\documentclass[11pt]{article}
\usepackage{amsthm,amsmath,amssymb}
\usepackage{indentfirst}
\usepackage{mathrsfs}
\usepackage{amsthm}
\usepackage{enumitem}  
\usepackage{mathrsfs}
\usepackage{latexsym,amssymb,amsmath}
\usepackage{graphicx,float,color,fancybox,shapepar,setspace,hyperref}
\usepackage{subfigure}
\usepackage{pgf,tikz}
\usepackage[affil-it]{authblk}
\usetikzlibrary{arrows}
\newtheorem{thm}{Theorem}[section]

\newtheorem{rk}[thm]{Remark}

\newtheorem{cor}[thm]{Corollary}
\newtheorem{lem}[thm]{Lemma}

\def\qed{\hfill\square}

\def\qed{ \hfill $\blacksquare$}

\newcounter{counter}
\newcommand{\Mod}[1]{\ (\mathrm{mod}\ #1)}

\begin{document}
	\title{An improvement on the parity of Schur’s partition function}

	\author{$^{1}$Yiwen Lu,\,$^{2}$Tao Wei,\,$^{3}$Xuejun Guo 	\thanks{The authors are supported by National Nature Science Foundation of China (Nos. 11971226, 12231009).}}
	\affil{ $^{1,2,3}$Department of Mathematics, Nanjing University, Nanjing 210093, China}
	\affil{$^{1}$luyw@smail.nju.edu.cn\ \ \  $^{2}$weitao@smail.nju.edu.cn\ \ \ $^{3}$guoxj@nju.edu.cn}
	
	\date{ }
	\maketitle

\begin{abstract}
We improve  S.-C. Chen's result on the parity of Schur’s partition function. Let $A(n)$ be the number of Schur’s partitions of $n$, i.e., the number of partitions of $n$ into
distinct parts congruent to $1,\; 2 \Mod{3}$. S.-C. Chen \cite{MR3959837} shows $\small \frac{x}{(\log{x})^{\frac{47}{48}}} \ll \sharp \{0\le n\le x:A(2n+1)\; \text{is odd}\}\ll \frac{x}{(\log{x})^{\frac{1}{2}}}$. In this paper, we improve Chen's result to
 $\frac{x}{(\log{x})^{\frac{11}{12}}} \ll \sharp \{0\le n\le x:A(2n+1)\; \text{is odd}\}\ll \frac{x}{(\log{x})^{\frac{1}{2}}}.$ 
 \\
 
\end{abstract}
\section{Introduction} 
Let $n \geqslant 1$ be an integer, a partition of $n$ is any non-increasing sequence of natural numbers whose sum is $n$. We denote by $p(n)$ the number of partitions for $n \geqslant 1$, and for convenience, let $p(0)=1$. 
\par
Obviously, it is impossible to obtain an exact expression for $p(n)$, and one seeks instead to figure out the asymptotic behavior of $p(n)$ as $n$ increases. The first breakthrough is attributed to G. H. Hardy and Ramanujan \cite{MR1575586}. Twenty years later, their results were refined by Hans Rademacher \cite{MR1575213} into the now well-known Hardy-Ramanujan-Rademacher Asymptotic Formula,
$$p(n)=\frac{1}{\pi \sqrt{2}} \sum_{k=1}^{\infty} A_{k}(n) \sqrt{k} \cdot \frac{d}{d n}\left(\frac{1}{\sqrt{n-\frac{1}{24}}} \sinh \left[\frac{\pi}{k} \sqrt{\frac{2}{3}\left(n-\frac{1}{24}\right)}\right]\right).$$
This formula can even be used to calculate the number of partitions for any specific integers $n \ge 1$. This result is so accurate that it was highly praised by Hardy.
We briefly introduce Rademacher’s result here. Rademacher constructed some explicit functions $T_q(n)$ satisfying $p(n)=\sum_{q=1}^{\infty } T_q(n)$ for all $n$ and then obtained the asymptotic formula $$p(n)\sim\frac{1}{4n\sqrt{3} } e^{\pi\sqrt{2n/3} },$$ by $T_1(n)$. Moreover, he got the exact error such that there are explicit constants $A$ and $B$ making $$\left | p(n)-\sum_{q=1}^{A\sqrt{n} } T_q(n) \right | <\frac{B}{n^{1/4}} .$$
\par
Even if the above work has been established, we have no evidence to believe that it possesses any attractive arithmetic properties. For example, there is no evidence that would lead us to believe that $p(n)$ should have a preference to be even rather than odd. Doing a computer calculation for the first 10,000 values, we get that there are 5004 odd values and 4996 even values, which means that the ratio is almost 1:1. Then replacing 2 with 3 we find that among the first 10,000 values, there are 3,313; 3,325; and 3,362 values that are congruent to 0, 1, and 2 modulo 3 respectively, in a ratio of roughly 1:1:1\cite{MR1854533}. O. Kolberg \cite{MR117213} has proved that $p(n)$ is infinitely often even and infinitely often odd while Parkin and Shanks \cite{MR227126} conjectured that
$$\sharp \{0\le n\le x:p(n)\; \text{is even (resp. odd)}\}\sim \frac{1}{2} x, \; \text{as} \; x\rightarrow  \infty .$$
\par
It could be seen that the above conjecture is extremely difficult to prove but there are still many special partition functions that are very attractive. The celebrated partition theorem which Schur \cite{MR111111} proved in 1926 is:
\begin{thm}[\cite{MR111111}]
Let $A(n)$ denote the number of partitions of $n$ into distinct parts $\equiv 1,\; 2 \Mod{3}$. Let $A_1(n)$ denote the number of partitions of $n$ with minimal difference $3$ between parts and such that no two consecutive multiples of $3$ occur
as parts. Then
$$A(n)= A_1(n)$$
\end{thm}
Then in 1971 Andrews \cite{MR222222} found the following companion to Schur's theorem by a computer search:
\begin{thm}[\cite{MR222222}]
Let $A_2(n)$ denote the number of partitions of $n$ in the form $n=e_1+e_2+\cdots +e_\nu $ such that $e_l-e_{l+1}\ge 3,\; 2 \; \text{or} \; 5$ if $e_l \equiv 1, \;2 \;\text{or} \;3 \Mod{3}$.
Then 
$$A(n) = A_1(n) = A_2(n).$$
\end{thm}
A result by S.-C. Chen \cite{MR3959837} tells us $\small \frac{x}{(\log{x})^{\frac{47}{48}}} \ll \sharp \{0\le n\le x:A(2n+1)\; \text{is odd}\}\ll \frac{x}{(\log{x})^{\frac{1}{2}}}.$ In this paper, we constructed a powerful theorem by Gauss' genus theory. Then by the constructed theorem, we finally covered and generalized Chen's result. The final result is as follows: 
\begin{thm}
We have
$$\frac{x}{(\log{x})^{\frac{11}{12}}} \ll \sharp \{0\le n\le x:A(2n+1)\; \text{is odd}\}\ll \frac{x}{(\log{x})^{\frac{1}{2}}}.$$
\end{thm}

\section{Main results} 
We let $R(n,ax^2+bxy+cy^2)$ be the number of the representations of $n$ by $ax^2+bxy+cy^2$ with $x,\; y \in \mathbb{Z}$. By SageMath, we find that the reduced primitive positive definite binary quadratic forms of discriminant $-216$ are
\begin{align*}
x^2+54y^2,\; 2x^2+27y^2,\; 5x^2\pm 2xy+11y^2,\; 7x^2\pm6xy+9y^2.
\end{align*}
In this paper, we consider square-free integers $m\equiv 11\Mod{24}$. Note that $m$ can not be represented by $x^2+54y^2$ or $7x^2 \pm 6xy+9y^2$ from the arithmetic of modulo 3. Therefore Dirichlet’s theorem on binary quadratic forms \cite[Theorem 1]{10.2307/2371469} shows that 
\begin{align*}
	R(m,2x^2+27y^2)+2R(m,5x^2+2xy+11y^2)=2\sum_{d \mid m} \left ( \frac{-216}{d}  \right ) = 2\sum_{d \mid m} \left ( \frac{-6}{d}  \right ) ,
\end{align*}
where $R(m,5x^2+2xy+11y^2)=R(m,5x^2-2xy+11y^2)$ follows the fact that a solution $(x_0,y_0)$ to $m = 5x^2+2xy+11y^2$ corresponds to a solution $(x_0,-y_0)$ to $m = 5x^2-2xy+11y^2$, and $\left ( \frac{\cdot}{\cdot}  \right )$ is the Jacobi–Kronecker symbol (for more details about this, see \cite{MR3959837}).
\par
Let $m=\prod_{i=1}^{t} p_i $  be the prime factorization of $m$. Then
\begin{align*}
	R(m,2x^2+27y^2)+2R(m,5x^2+2xy+11y^2) = 2\sum_{d \mid m} \left ( \frac{-6}{d}  \right ) =\prod_{i=1}^{t} \left ( 1+\left ( \frac{-6}{p_i}  \right )  \right ) .
\end{align*}
It follows immediately that an arbitrary prime factor $p_i$ of $m$ must satisfy $\left ( \frac{-6}{p_i}  \right ) = 1$ if $ R(m,5x^2+2xy+11y^2)\ne 0 $, hence 
\begin{align}\label{formula2}
	p_i\equiv 1, 5, 7 ,11\Mod{24},\quad 1\leq i \leq t,
\end{align}
and we get 
\begin{align}\label{formula3}
	R(m,2x^2+27y^2)+2R(m,5x^2+2xy+11y^2) =2^{t+1}.
\end{align}
A similar discussion for the discriminant $-24$ (see also \cite{MR3959837}) tells us that
\begin{align}\label{formula77}
	R(n,2x^2+3y^2)+R(n,x^2+6y^2)=2\sum_{d \mid n} \left ( \frac{-6}{d}  \right ), \quad n\in \mathbb{N},
\end{align}
and by considering module 3, we have
\begin{align}\label{formula101}
R(m,2x^2+3y^2)=2^{t+1}.
\end{align}
\par
We denote by $C(D)$ the class group (see \cite[p.\,45-46]{MR3236783} for definitions) of discriminant $D$ and let $h(D)$ denote the number of classes of primitive positive definite forms of discriminant $D$. By \cite[Table 9.1]{MR2234833}, we know $h(-24)=2$ and $2x^2+3y^2$ is the generator in $C(-24)$.\\
For convenience, we write
\begin{align*}
	f(x,y)=x^2+6y^2, \quad g(x,y)=2x^2+3y^2.
\end{align*}
Thus we get the Dirichlet compositions (see \cite{MR3236783} for definitions) as follows:
\begin{equation}\label{equation}
	\begin{aligned}
		f(x,y)f(z,\omega)&=f(xz-6y\omega,x\omega+yz)\\
		f(x,y)g(z,\omega)&=g(xz-3y\omega,x\omega+2yz)\\
		g(x,y)g(z,\omega)&=f(2xz-3y\omega,x\omega+yz).
	\end{aligned}
\end{equation}
\par
In this paper, we used Gauss' genus theory to construct a powerful theorem.  Our main results are as follows:
\begin{thm}\label{th1}
Let  $m\equiv11\Mod{24}$ be a square-free integer. Then \begin{align*}R(m,5x^2+2xy+11y^2) \equiv 2\Mod{4},
	\end{align*}
if and only if $m$ can be written in the following two forms:
\begin{enumerate}[label=(\roman*)]
\item 
$m=p_1 \cdots p_{2t_1-1} \cdot \overline{p}_1 \cdots \overline{p}_{2t_2} \cdot q_1 \cdots q_{t_3} \cdot \overline{q}_1 \cdots \overline{q}_{2t_4};$
\item 
$m=p_1 \cdots p_{2t_1} \cdot \overline{p}_1 \cdots \overline{p}_{2t_2-1} \cdot q_1 \cdots q_{t_3} \cdot \overline{q}_1 \cdots \overline{q}_{2t_4-1},$
\end{enumerate}
where $p_1, \cdots ,p_{2t_1} \in \mathscr{S}_1$, $\overline{p}_1, \cdots ,\overline{p}_{2t_2} \in \mathscr{S}_2$, $q_1, \cdots ,q_{t_3} \in \mathscr{S}_3$, $\overline{q}_1, \cdots ,\overline{q}_{2t_4} \in \mathscr{S}_4$, $t_1,\;t_2,\;t_3,\;t_4 \in \mathbb{N}$, and $\mathscr{S}_1$, $\mathscr{S}_2$, $\mathscr{S}_3$, $\mathscr{S}_4$ are the following four subsets of primes:
\begin{align*}
\mathscr{S}_1&=\{p:p\equiv 11\Mod{24}, p=5x^2+2xy+11y^2\}
\\
 \mathscr{S}_2&=\{\overline{p}:\overline{p}\equiv 5\Mod{24}, \overline{p}=5x^2+2xy+11y^2\}
\\
\mathscr{S}_3&=\{q:q\equiv 1\Mod{24}, q=7x^2+6xy+9y^2\}
\\
 \mathscr{S}_4&=\{\overline{q}:\overline{q}\equiv 7\Mod{24}, \overline{q}=7x^2+6xy+9y^2\}.
\end{align*}
 \end{thm}\medskip
Then by Theorem \ref{th1}, we covered and generalized the result of S.-C. Chen \cite{MR3959837}. The contents are as follows:
\begin{thm}\label{th2}
	We have
$$\frac{x}{(\log{x})^{\frac{11}{12}}} \ll \sharp \{0\le n\le x:A(2n+1)\; \text{is odd}\}\ll \frac{x}{(\log{x})^{\frac{1}{2}}}.$$
	\end{thm}
Since $\sharp\{1\le n\le x:A(n)\;\text{is odd}\}\ge \sharp\{0\le n\le\frac{x-1}{2} :A(2n+1)\;\text{is odd}\}$, we have the following corollary.
  \begin{cor}
 $$\sharp\{1\le n\le x:A(n)\;\text{is odd}\}\gg \frac{x}{(\log{x})^{\frac{11}{12}}} .$$
  \end{cor}
\section{Proofs}
Before starting our proofs of the main theorems, there is a lemma which will be used in the proof of Theorem \ref{th1}.
\begin{lem}\label{lem} 
	For the Dirichlet compositions in \eqref{equation} and any square-free integers $a=f(x_0,y_0)\ne 1$ or $a=g(x_0,y_0)\ne 1$, $b=f(x_1,y_1)\ne 1$ or $g(x_1,y_1) \ne 1$
	where $3 \nmid x_i$ and $x_i, \;y_i\in \mathbb{N},\; i=0,1.$ If
	$$ab=f(x_2,y_2)\;(resp.\;g(x_2,y_2))\; \text{and} \; ab=f(x_2^{\prime},y_2^{\prime})\;(resp.\;g(x_2^{\prime},y_2^{\prime})),$$
	where $(x_2,y_2)$ and $(x_2^{\prime},y_2^{\prime})$ are composed of $(x_0,y_0)$ and $(x_1,y_1)$,
	then
	\begin{enumerate}[label=(\roman*)]
		\item 
		$3\mid y_0,\; 3\mid y_1 \Leftrightarrow 3 \mid y_2,\;3 \mid y_2^{\prime}$;
		\item 
		$3\nmid y_0,\; 3\nmid y_1 \Leftrightarrow 3\mid y_2,\; 3\nmid y_2^{\prime}$ or $3\nmid y_2,\; 3\mid y_2^{\prime};$
		\item 
		$3\mid y_0,\; 3\nmid y_1$ or $3\nmid y_0,\; 3\mid y_1 \Leftrightarrow 3 \nmid y_2,\;3 \nmid y_2^{\prime}.$
	\end{enumerate}
\end{lem}
\noindent{\bf Proof of Lemma \ref{lem}.} We take the Dirichlet composition of $f(x,y)$ and $f(x,y)$ as an example first. By equation \eqref{equation}, we know
$$y_2=x_0y_1+y_0x_1,\; y_2^{\prime}=x_0y_1-y_0x_1\;(resp.\;y_2=x_0y_1-y_0x_1,\; y_2^{\prime}=x_0y_1+y_0x_1) .$$
The necessity of (i) is obvious. Conversely, we get that $3\mid 2x_0y_1 = y_2+y_2^{\prime}$ and $3 \mid 2y_0x_1=\left |y_2-y_2^{\prime}\right |$. Thus $ 3\mid y_i ;$ Note that $y_2y_2^{\prime} = x_{0}^2y_{1}^2-y_{0}^2x_{1}^2$. Combing the fact $x_{i}^2 \equiv 1 \Mod{3}$, $y_{i}^2 \equiv 1 \Mod{3}$ and (i), the necessity of (ii) is proved. Conversely, we have $3 \mid y_2y_2^{\prime} =x_{0}^2y_{1}^2-y_{0}^2x_{1}^2$, i.e., $y_{1}^2-y_{0}^2 \equiv 0 \Mod{3}$. By (i) again, the sufficiency is proved; The necessity of (iii) is proved immediately according to the fact that $3\nmid y_{1}^2-y_{0}^2$. Then combining (i) and (ii), the converse holds naturally.
\par
Note that $x_0y_1+2y_0x_1 \equiv x_0y_1-y_0x_1 \Mod{3}$ and $x_0y_1-2y_0y_1 \equiv x_0y_1+y_0x_1 \Mod{3}$. Then a similar discussion tells us that the same conclusion holds in the Dirichlet composition of $f(x,y)$ and $g(x,y)$ (resp. $g(x,y)$ and $g(x,y)$).\qed\\

\noindent {\bf Proof of Theorem \ref{th1}.} For the case that $m$ is a prime, we have $R(m,2x^2+27y^2)+2R(m,5x^2+2xy+11y^2)=4$ by \eqref{formula3}. Obviously
\begin{align}\label{formula33}
 R(m,5x^2+2xy+11y^2)\equiv 2 \Mod{4} \Leftrightarrow 	R(m,5x^2+2xy+11y^2)>0 ,
\end{align}
i.e., $m \in \mathscr{S}_1$ and can be written as the first form in Theorem \ref{th1}.
\par
For the case that $m$ is not a prime, we deduce from \eqref{formula3} that:
$$\frac{1}{4} R(m,2x^2+27y^2)\equiv \frac{1}{2}R(m,5x^2+2xy+11y^2)\Mod{2},$$
which states that
\begin{align}\label{formula34}
R(m,5x^2+2xy+11y^2)\equiv 2 \Mod{4} \Leftrightarrow R(m,2x^2+27y^2)=4k,\; \text{for an odd integer}\, k .
\end{align}
\noindent {\bf Claim\refstepcounter{counter} \arabic{counter}.}
$R(m,2x^2+27y^2)=4k$, for an odd integer $k$ in this case if and only if $m$ is the form of (i) or (ii) in Theorem \ref{th1}. 
\par
We find that a solution $(a,b)$ to $m=2x^2+27y^2$ must be a solution $(a,3b)$ to $m=2x^2+3y^2$. Hence 
\begin{align}\label{three}
	\sharp\{(a,b):2a^2+27b^2=m,a\in \mathbb{N},b\in \mathbb{N}\}=\sharp\{(a^{\prime},b^{\prime}):2a^{\prime2}+3b^{\prime2}=m,a^{\prime}\in \mathbb{N},b^{\prime}\in \mathbb{N}, 3\mid b^{\prime}\},
\end{align}
which means we can explore the solutions of $m=2x^2+3y^2$ that the corresponding $y$ value is divisible by $3$ for proving the claim. 
\par
Combing the fact that a prime $p \equiv 5, \;11 \Mod{24}$ can be represented by $2x^2+3y^2$, a prime $q \equiv 1,\; 7 \Mod{24}$ can be represented by $x^2+6y^2$ \cite[Table 9.1]{MR2234833} and formula \eqref{formula77}, we get
\begin{align}
		&R(p_i,2x^2+3y^2)=4, \; \text{if} \; p_i \equiv 5,\; 11\Mod{24}, \nonumber \\
	&R(p_i,x^2+6y^2)=4, \;\text{if} \;  p_i \equiv 1,\; 7\Mod{24} \nonumber.
\end{align}
Therefore after adjusting the order of the prime factor of $m$, $m$ can be written as
\begin{equation}
\begin{aligned}\label{equation m}
m=\prod_{i=1}^{i_0} (x_1^2 + 3 y_1 ^2)\cdots(2  x_{i_0}^2+3   y_{i_0}^2)  \cdot\prod_{i=i_0+1}^{t}(  x_{i_0+1}^2+6  y_{i_0+1}^2)\cdots( x_t^2+6  y_t^2),
\end{aligned}
\end{equation}
where 
\begin{eqnarray}\label{mi}
	p_i=\left\{
	\begin{array}{rcl}
		2x_i^2+3y_i^2 & & {1 \leq i \leq i_0},\\
		x_i^2+6y_i^2 & & {i_0 < i \leq t},
	\end{array}
\right.
\end{eqnarray}
where $i_0 =
\#\{p_i:p_i \equiv 5,\; 11 \Mod{24}\}$.
\par
It is easy to see that an arbitrary solution $(x_0,y_0)$ to equation $m=2x^2+3y^2$ corresponds to three other solutions $(x_0,-y_0)$, $(-x_0,-y_0)$ and $(x_0,-y_0)$. We call such four solutions a class of solutions to the equation $m=2x^2+3y^2$. Therefore the equation $m=2x^2+3y^2$ has $2^t-1$ classes of solutions by \eqref{formula101}. Similarly, $2x^2+3y^2=p_i \equiv 5, 11\Mod{24}$ and 	$x^2+6y^2=p_i\equiv 1, 7\Mod{24}$ both have a class of solutions. Then by \eqref{equation}, we know that an arbitrary class of solutions to $f(x,y)=a$ and an arbitrary class of solutions to $f(x,y)=b$ where
$$a \mid m,\; b \mid m,\;ab \mid m,\; a\ne1,\; b\ne1 ,$$
correspond to two classes of solutions to $f(x,y)=ab$ as follows:
$$
	[(\pm  x_0 )^2+6(\pm  y_0 )^2 ]\cdot[(\pm  x_1 )^2+6(\pm  y_1 )^2] =[\pm  (x_0x_1-6y_0y_1 )]^2+6[\pm (x_0y_1+y_0x_1 )]^2,
$$
$$
	[(\pm  x_0 )^2+6(\pm  y_0 )^2 ]\cdot[(\pm  x_1 )^2+6(\pm  y_1 )^2] =[\pm ( x_0x_1+6y_0y_1) ]^2+6[\pm ( x_0y_1-y_0x_1)]^2.
$$
Similarly, we have:
$$[(\pm  x_0 )^2+6(\pm  y_0 )^2 ]\cdot[2(\pm  x_1 )^2+3(\pm  y_1 )^2] =2[\pm ( x_0x_1-3y_0y_1)]^2+3[\pm (x_0y_1+2y_0x_1 )]^2,$$
$$[(\pm  x_0 )^2+6(\pm  y_0 )^2 ]\cdot[2(\pm  x_1 )^2+3(\pm y_1 )^2] =2[\pm ( x_0x_1+3y_0y_1 )]^2+3[\pm ( x_0y_1-2y_0x_1 )]^2,$$
which corresponds to the Dirichlet composition of $f(x,y)$ and $g(x,y)$,
\\
and
$$[2(\pm  x_0 )^2+3(\pm  y_0 )^2][2(\pm  x_1 )^2+3(\pm  y_1 )^2]=[\pm ( 2x_0x_1-3y_0y_1)]^2+6[\pm (x_0y_1+y_0x_1 )]^2,$$
$$[2(\pm  x_0 )^2+3(\pm  y_0 )^2][2(\pm  x_1 )^2+3(\pm  y_1 )^2]=[\pm ( 2x_0x_1+3y_0y_1 )]^2+6[\pm (x_0y_1-y_0x_1 )]^2,$$
which corresponds to the Dirichlet composition of $g(x,y)$ and $g(x,y)$. 
\par
We call the corresponding discussed above the composition of solutions in the process of Dirichlet compositions. The above discussion states that an arbitrary solution to the equation $m=2x^2+3y^2$ can be obtained via the composition of solutions that in the process of Dirichlet compositions in formula \eqref{equation} step by step. 
\\
\noindent {\bf Claim\refstepcounter{counter}\label{smallpart}  \arabic{counter}.}\label{claim2}
Now we claim that only all pairs $(x_i ,y_i)$ extracted from the right side of equation \eqref{equation m} satisfies $3\nmid y_i$, then $R(m,2x^2+27y^2)=4k$, for an odd integer $k$.
\par
We assume that there are $t_0$ ($t_0 \in \mathbb{N}$) pairs $(x_i ,y_i)$ in equation \eqref{equation m} that $3 \mid y_i$ while the rest of $t-t_0$ pairs $(x_i ,y_i)$ satisfies $3 \nmid y_i $.  Then we adjust the order of $p
_i$ to make the $4t_o$ pairs $(\pm x_i ,\pm y_i)$ in the right side of equation \eqref{equation m} conduct the composition of solutions firstly. By lemma \ref{lem}, we finally get $2^{t_0 +1}$ pairs $(x,y)$ $(x,\; y \in \mathbb{Z})$ that satisfies
$y \equiv 0 \Mod{3}$ and 
$$2x^2+3y^2=p_1\cdots p_{t_0}, \quad \text{or} \quad x^2+6y^2=p_1\cdots p_{t_0},$$
where $p_1,\cdots,p_{i_0}$ are the corresponding primes of $t_0$ pairs $(x_i ,y_i)$ above.
\par
Then we let the right side of above equation be multiplied one by one with the remaining $t-t_0$ $p_i$ (that is, let the $2^{t_0+1}$ pairs we have obtained with the remaining  $4(t-t_0)$ pairs $(\pm x_i, \pm y_i)$ conduct the composition of solutions). We denote by $u(n)$ the number of the pairs $(x,y)$ $(x,\; y \in \mathbb{Z})$
satisfying $3 \mid y$ and denote by $v(n)$ the number of the pairs $(x,y)$ $(x, \;y \in \mathbb{Z})$ satisfying $3 \nmid y$ after conducting the composition of solutions $n$ $(1\leq n \leq t-t_0)$ times with the remaining $t-t_0$ $p_i$ in turn. Therefore
$$u(0)=2^{t_0+1}, \quad v(0)=0,$$
and
\begin{align*}
	u(n)=v(n-1), \quad v(n)=2u(n-1)+v(n-1), \quad 1\leq n \leq t-t_0,
\end{align*}
which follows from lemma \ref{lem}. We immediately get that
\begin{equation}\label{equation8}
\begin{aligned}
	&u(t-t_0)=v(t-t_0-1)=2u(t-t_0-2)+v(t-t_0-2)\\
	&=2u(t-t_0-3)+2u(t-t_0-3)+v(t-t_0-3)\\
	&=\cdots\\
	&=2(c_1 u(1)+c_2 v(1))+v(1)\\
	&=c_2 2^{t_0+3}+2^{t_0+1}
\end{aligned}
\end{equation}
where $c_1, \; c_2 \in \mathbb{N}$.\\
Putting \eqref{three} and \eqref{equation8} together, we obtain
$$t_0=0 \Leftrightarrow R(m,2x^2+27y^2)=4k,\; k \;\text{is an odd integer} .$$
Hence Claim2 is proved and we will prove Claim1 by Claim2 next.
\par
Now by Claim2, we know that for all prime factors $p_i$ of $m$,
$$R(m,2x^2+27y^2)=4k,\; k \;\text{is an odd integer},$$
if and only if for all $y_i$ in formula \eqref{mi}, $3\nmid y_i$. This shows that 
\begin{equation}\label{mii}
	\begin{aligned}
			&R(p_i,2x^2+27y^2)=0, \; \text{if} \; p_i \equiv 5, \;11\Mod{24},\\
		&R(p_i,x^2+54y^2)=0, \;\text{if} \;  p_i \equiv 1,\; 7\Mod{24}, 
	\end{aligned}
\end{equation}
for all prime factors $p_i$ of $m$. Note that $p_i \equiv 5,\; 11\Mod{24}$ can not be represented by $x^2+54y^2$ or $7x^2 \pm 6xy+9y^2$ by considering modulo 3. Therefore combining Dirichlet’s theorem on binary quadratic forms \cite[Theorem 1]{10.2307/2371469} and \eqref{mii}, we have 
$$R(p_i,5x^2+2xy+11y^2)=2, \; \text{if} \; p_i \equiv 5,\; 11\Mod{24}.$$
Note that $p_i \equiv 1, \;7\Mod{24}$ can not be represented by $2x^2+27y^2$ and $5x^2 \pm 2xy+11y^2$ by considering modulo 3, similarly we have
$$R(p_i,7x^2+6xy+9y^2)=2, \; \text{if} \; p_i \equiv 1,\; 7\Mod{24}.$$
Thus 
$$\forall p_i\mid m,  \; p_i \in \mathscr{S}_1\cup \mathscr{S}_2 \cup\mathscr{S}_3 \cup \mathscr{S}_3.$$
Moreover, the fact $m \equiv 11 \Mod{24}$ ensures that $m$ can only be of the two forms in Theorem \ref{th1}.
Claim1 is proved.
\par
If $R(m,5x^2+2xy+11y^2) \equiv 2 \Mod{4}$, we have \eqref{formula2} and equation \eqref{formula3}. Hence by \eqref{formula33}, \eqref{formula34} and Claim1, $m$ must be the form of (i) or (ii) in Theorem \ref{th1}. Conversely, equation \eqref{formula3} also holds. Putting  \eqref{formula33}, \eqref{formula34} and Claim1 together, we obtain $R(m,5x^2+2xy+11y^2) \equiv 2 \Mod{4}$. 
This completes the proof of Theorem~\ref{th1}.\qed\\
\\
\noindent{\bf Proof of Theorem \ref{th2}.} Define $\mathscr{S}_1$, $\mathscr{S}_2$, $\mathscr{S}_3$ and $\mathscr{S}_4$ as the subsets of primes in Theorem \ref{th1}, respectively. We let $\mathscr{S}=\mathscr{S}_1 \cup \mathscr{S}_2 \cup \mathscr{S}_3 \cup \mathscr{S}_4$ and 
$$\mathscr{C}=\{x:x \in \mathbb{N}, \,x\;\text{is square-free },\, p\mid x \Rightarrow p \in \mathscr{S}\}.$$
For simplicity, we define several functions next. If $n$ is a square-free integer, then we let $\mu_i(n)$=
\begin{eqnarray}\nonumber
	\left\{
	\begin{array}{rcl}
	1 & & {\text{if the number of factors from} \,\mathscr{S}_i \,\text{in the prime factorization of} \,n\,\text{is even} };\\
		-1 & & {\text{if the number of factors from} \,\mathscr{S}_i \,\text{in the prime factorization of} \,n\,\text{is odd}},
	\end{array}
	\right.
\end{eqnarray}
where $1 \leqslant i \leqslant 4.$ By \cite[(8)]{MR3959837}, we have
\begin{align}\label{chen}
	A(2n+1)\equiv \frac{1}{2} R(24n+1,5x^2+2xy+11y^2)\Mod{2}.
\end{align}
Combing \eqref{chen} with Theorem \ref{th1}, we get 
$$A\left ( \frac{m+1}{12}  \right ) \equiv 1 \Mod{2},$$
for any square-free integer $m$ that is the form of (i) or (ii) in Theorem \ref{th1}. Thus
\begin{equation}\label{inequality14}
\begin{aligned}
	\sum_{\substack{0\le n\le x\\A(2n+1)\; \text{odd}}}1&\ge 	\sum_{\substack{m\le x \\ m \;\text{is the form (i)}}}1 + \sum_{\substack{m\le x \\ m \;\text{is the form (ii)}}}1\\
	&=\sum_{\substack{m\in \mathscr{C},\,m\le x \\ \mu_1 (m)=-1,\mu_2 (m)=1,\mu_4 (m)=1}}1 
		+
		\sum_{\substack{m\in \mathscr{C},\,m\le x \\ \mu_1 (m)=1,\mu_2 (m)=-1,\mu_4 (m)=-1}}1 
\end{aligned}
\end{equation}
Since the number of classes of discriminant $-216$ is $6$, the Chebotarev density theorem \cite[Theorem 9.12]{MR3236783} shows that the Dirichlet density of the set of primes represented by $5x^2+2xy+11y^2$ is $\frac{1}{6}$. Applying the orthogonality of Dirichlet character modulo $24$, we see that the Dirichlet density of $\mathscr{S}_1$ is $\frac{1}{6} \cdot\frac{1}{\phi (24)} =\frac{1}{48}$, where $\phi $ is Euler’s totient function. Similarly,  the Dirichlet density of $\mathscr{S}_2$, $\mathscr{S}_3$ and $\mathscr{S}_4$ are all $\frac{1}{6} \cdot\frac{1}{\phi (24)} =\frac{1}{48}$. Thus the Dirichlet density of $\mathscr{S}$ is $4\cdot\frac{1}{48} =\frac{1}{12}$.
\\
We define
$$h(x)=\sum_{m\in \mathscr{C},\,m\le x }1,$$
and
\begin{align*}
	 h_{1} (x)=\sum_{\substack{m\in \mathscr{C},\,m\le x \\ \mu_1 (m)=1}}1 , \qquad 
		h_{-1} (x)=\sum_{\substack{m\in \mathscr{C},\,m\le x \\ \mu_1 (m)=-1}}1 ,
\end{align*}
where $h(x)=h_{1} (x)+h_{-1} (x)$ obviously. A classical result of Wirsing \cite{MR333333} on multiplicative functions (see also \cite[Proposition 4]{MR2644888}) tells us that 
\begin{align}\label{formula16}
	h(x)=C_h \frac{x}{(\log{x})^\frac{11}{12}} + o\left (\frac{x}{(\log{x})^\frac{11}{12}}\right ),
\end{align}
where constant $C_h>0$, as $x\to \infty$.
\\
Note that $h_{1}(x) \geq 0$,  $h_{-1}(x) \geq 0$. We deduce from \eqref{formula16} that
$$h_{1}(x)=C_{h_1} \frac{x}{(\log{x})^\frac{11}{12}} + o \left (\frac{x}{(\log{x})^\frac{11}{12}} \right ),\qquad h_{-1}(x)=C_{h_{-1}} \frac{x}{(\log{x})^\frac{11}{12}} + o\left (\frac{x}{(\log{x})^\frac{11}{12}}\right) ,$$
where there are three possible cases of $C_{h_1}$ and $C_{h_{-1}}$ as follows:
\begin{enumerate}[label=(\arabic*)]
	\item 
	$C_{h_1}>0$, $C_{h_{-1}}>0$ and $C_{h_1} + C_{h_{-1}} = C_h$;
	\item 
	$C_{h_1}=C_h$, $C_{h_{-1}}=0$;
	\item 
	 $C_{h_1}=0$, $C_{h_{-1}}=	C_h$.
\end{enumerate}
We find that
\begin{equation}
	\begin{aligned}\label{formula17}
		\small h_{1} (x)&=
		\small \sum_{\substack{m\in \mathscr{C},\,m\le x \\ 11\mid m,\,\mu_1 (m)=1}}1 \small
		+
		\small \sum_{\substack{m\in \mathscr{C},\,m\le x \\ 11\nmid m,\, \mu_1 (m)=1}}1 \small
		\\
		&=\small \sum_{\substack{m\in \mathscr{C},\,m\le \frac{x}{11} \\ 11\nmid m,\, \mu_1 (m)=-1}}1 \small
		+ \small
    	\sum_{\substack{ m\in \mathscr{C},\,m\le 11x \\ 11\mid m,\, \mu_1 (m)=-1}}1 \small	\\
		&\le \small \sum_{\substack{m\in \mathscr{C},\,m\le \frac{x}{11} \\  \mu_1 (m)=-1}}1 \small
		+ \small
		\sum_{\substack{m\in \mathscr{C},\,m\le 11x \\  \mu_1 (m)=-1}}1 \small\\
		&=h_{-1} (\frac{x}{11})+h_{-1} (11x).						
	\end{aligned}
\end{equation}
Similarly, we have
\begin{align}\label{formula18}
h_{-1}(x) \le h_{1}(\frac{x}{11})+h_{1} (11x).
\end{align}
Therefore case (1) holds because case (2) and case (3) are absurd to the inequality \eqref{formula17} and inequality \eqref{formula18}, respectively, i.e.,
$$\sum_{\substack{m\in \mathscr{C},\,m\le x \\ \mu_1 (m)=-1}}1= h_{-1} (x)=C_{h_{-1}} \frac{x}{(\log{x})^\frac{11}{12}} + o \left (\frac{x}{(\log{x})^\frac{11}{12}} \right ) ,$$
where constant $C_{h_{-1}}>0$, as $x\to \infty$.
Then we write $h_1(x)$ above still as $h(x)$. In a similar way, we construct the new functions $h_{1} (x)$ and $h_{-1} (x)$ as follows:
\begin{align*}
	\small h_{1} (x)=\sum_{\substack{m\in \mathscr{C},\,m\le x \\ \mu_1 (m)=-1,\,\mu_2 (m)=1}}1,\qquad 
	\small h_{-1} (x)=\sum_{\substack{m\in \mathscr{C},\,m\le x \\ \mu_1 (m)=-1,\,\mu_2 (m)=-1}}1,
\end{align*}
where $h_{1} (x) + h_{-1} (x)= h(x)$ and 
$$h_{1}(x) \le h_{-1}(\frac{x}{5})+h_{-1} (5x),$$
$$h_{-1}(x) \le h_{1}(\frac{x}{5})+h_{1} (5x).$$
Similarly, we get 
\begin{align*}
	\small
\sum_{\substack{m\in \mathscr{C},\,m\le x \\ \mu_1 (m)=-1,\,\mu_2 (m)=1}}1=h_{1}(x)=C_0\frac{x}{(\log{x})^\frac{11}{12}} + o \left (\frac{x}{(\log{x})^\frac{11}{12}} \right ),
	\end{align*}
where constant $C_0>0$, as $x\to \infty$.
Continuously through similar discussions, we can get 
$$\small\sum_{\substack{m\in \mathscr{C},\,m\le x \\ \mu_1 (m)=-1,\mu_2 (m)=1,\mu_4 (m)=1}}1 \small = C_1\frac{x}{(\log{x})^\frac{11}{12}} + o \left (\frac{x}{(\log{x})^\frac{11}{12}} \right ),$$ and 
$$\small\sum_{\substack{m\in \mathscr{C},\,m\le x \\ \mu_1 (m)=1,\mu_2 (m)=-1,\mu_4 (m)=-1}}1 \small = C_2\frac{x}{(\log{x})^\frac{11}{12}} + o \left (\frac{x}{(\log{x})^\frac{11}{12}} \right ),$$
where $C_1>0$, $C_2>0$ as $x\to \infty$.
Then by \eqref{inequality14} we finally get
$$\sum_{\substack{0\le n\le x\\A(2n+1)\; \text{odd}}}1 \gg  \frac{x}{(\log{x})^\frac{11}{12}}.$$
Combining the the upper bound \cite[Thorem 1.1]{MR3959837} with the above inequality, we complete the proof of Theorem~\ref{th2}.\qed\\
\begin{rk}
Obviously, Theorem \ref{th1} can also be used to discuss the case that $m$ is not a square-free integer, i.e., the primes that appear in form (i) and (ii) in Theorem \ref{th1} can be the same. This allows one to obtain more integers $n \in \mathbb{N}$ that satisfy $ A(2n+1)\; \text{is odd} $. Therefore one may be able to further raise the lower bound of $\sharp \{0\le n\le x:A(2n+1)\; \text{is odd }\}$
	and obtain a more accurate asymptotic.
\end{rk}

\bibliographystyle{siam}

\end{document}